\documentclass{amsart}

\usepackage{mathtools}
\usepackage{fontspec}
\usepackage{unicode-math}
\usepackage{fullpage}

\usepackage{enumitem}

\usepackage{tikz-cd}
\tikzcdset{
  arrow style=tikz
}

\usepackage[dvipsnames]{xcolor}

\usepackage{comment}
\usepackage[style=american]{csquotes}

\usepackage[
  backend=biber,
  style=alphabetic,
  sorting=nyt,
  minnames=4,
  maxnames=5,
  mincitenames=4,
  maxcitenames=5,
  minalphanames=4,
  maxalphanames=5
]{biblatex}
\AtBeginBibliography{%
  \setlength{\emergencystretch}{3em}%
}

\usepackage{hyperref}
\hypersetup{
  pdfdisplaydoctitle = true,
  pdftitle           = { Autoformalizing the calculation of
                         \texorpdfstring{$\pi_3(S^2)$}{π₃(S²)}
                       },
  pdfauthor          = { Daniel Caranza and
                         Chunyi Liu and
                         Emily Riehl and
                         Egbert Rijke
                       }
}

\DeclareFieldFormat{url}{\allowbreak\url{#1}}

\usepackage{zref-clever}

\zcsetup{
  cap,
  noabbrev,
  nameinlink,
}

\zcRefTypeSetup{section}{
  Name-sg = {\S},
  name-sg = {\S},
  Name-pl = {\S\S},
  name-pl = {\S\S},
  namesep = {},
}

\theoremstyle{definition}

\theoremstyle{remark}

\makeatletter
\let\c@equation\c@thm
\makeatother
\numberwithin{equation}{section}

\usepackage{marginnote}

\newcommand{\ZZ}{\mathbb{Z}}

\newcommand{\lean}{Lean}
\newcommand{\mathlib}{Mathlib}
\newcommand{\leaneval}{lean-eval}
\newcommand{\agda}{Agda}
\newcommand{\agdaunimath}{agda-unimath}
\newcommand{\codex}{Codex}
\newcommand{\claude}{Claude}
\newcommand{\automath}{Automath}
\newcommand{\mizar}{Mizar}
\newcommand{\rocq}{Rocq}
\newcommand{\coq}{Coq}
\newcommand{\isabelle}{Isabelle}
\newcommand{\hollight}{HOL Light}

\newcommand{\supst}{\textsuperscript{st}}

\newcommand{\suprd}{\textsuperscript{rd}}
\newcommand{\supth}{\textsuperscript{th}}

\newcommand{\join}{\mathbin{\ast}}

\newcommand{\pt}{\mathsf{pt}}
\newcommand{\loopcircle}{\mathsf{loop}}
\newcommand{\north}{\mathsf{n}}
\newcommand{\south}{\mathsf{s}}
\newcommand{\merid}{\mathsf{merid}}
\newcommand{\ap}{\mathsf{ap}}
\newcommand{\refl}{\mathsf{refl}}
\newcommand{\unit}{\mathbf{1}}
\newcommand{\Fin}[1]{\mathsf{Fin}_{#1}}

\tikzset{alt/.code={}}

\usepackage{afterpage}
\usepackage{floatpag}

\usepackage{fancyvrb}
\usepackage{newunicodechar}
\newfontface\AgdaArxivMono{JuliaMono-Regular.ttf}[Scale=MatchLowercase]

\newunicodechar{＝}{=}
\newunicodechar{𝕊}{\ensuremath{\mathbb{S}}}

\DefineVerbatimEnvironment{agdacode}{Verbatim}{
  formatcom=\AgdaArxivMono,
  fontsize=\small
}
\CustomVerbatimCommand{\agdainline}{Verb}{
  formatcom=\AgdaArxivMono
}

\newcommand{\AgdaCodePlain}{}

\newcommand{\AgdaCodeName}[2]{}
\newcommand{\AgdaCodeSetColor}[2]{}
\newcommand{\AgdaCodeDisplayAs}[2]{}
\newcommand{\AgdaCodeRemoveName}[1]{}

\newcommand{\AgdaCodeSetup}[1]{}

\begin{document}

\title{Autoformalizing the calculation of \texorpdfstring{$\pi_3(S^2)$}{π₃(S²)}}

\author{Daniel Carranza}
\author{Chunyi Liu}
\author{Emily Riehl}
\author{Egbert Rijke}

\date{\today}

\address{Johns Hopkins University \\ 3400 N Charles Street \\ Baltimore, MD 21218}
\email{dcarran3@jhu.edu}
\email{cliu238@jhu.edu}
\email{eriehl@jhu.edu}
\email{erijke1@jhu.edu}

\begin{abstract}
  We report our findings and the results of our experiment of autoformalizing the homotopy type theoretic computation of $\pi_3(S^2)$ in \agdaunimath{} using \codex{}.
\end{abstract}

\maketitle


\section{Introduction}

In April 2026, the \lean{} Focused Research Organization released a \lean{} AI formalization benchmark \leaneval{} \cite{lean-eval}. The purpose of \leaneval{} is to evaluate AI-powered autoformalization systems that are designed to attempt computer formalizations of challenging mathematical problems. The benchmark problems were chosen so as to have \lean{}-formalizable theorem statements using \mathlib{}, \lean{}'s mathematics library \cite{mathlib}, with known natural language proofs that would take significant work to formalize.

One of the benchmark problems was to demonstrate an isomorphism $\pi_3(S^2) \cong \ZZ$, calculating the third homotopy group of the $2$-sphere \cite{lean-eval-pi3-s2}. This calculation involves some beautiful mathematical ideas: the generating element in the group $\ZZ$ corresponds under the isomorphism to a continuous map $S^3 \to S^2$ known as the \emph{Hopf fibration}, a fiber bundle with fiber $S^1$ (see \cite{johnson-hopf} for a visualization). This problem attracted our interest as a potential test of autoformalization capabilities in the proof assistant \agda{} using the \agdaunimath{} library, a library of mathematical statements proven in homotopy type theory \cite{agda-unimath}.

While \lean{} is the most popular interactive theorem prover used by mathematicians today, it is not the only option, and historic successes have been achieved with many proof assistants. Perhaps the first great achievement in formalized mathematics was the verification of Landau's \emph{Grundlagen} in \automath{} \cite{Benthem1979Checking}. Other notable achievements include the formalization of the Jordan curve theorem in \mizar{} \cite{Korniowicz:jordan}, the formalizations of the four-color theorem and the Feit--Thompson theorem in \rocq{} (\emph{nee}.~\coq{}) \cite{gonthier:feit-thompson}, and the formalization of the Kepler conjecture in \hollight{} and \isabelle{} \cite{hales:kepler}.

An interesting possibility enabled by the computer proof assistant paradigm, described in \cite{Shulman2024Strange}, is to write proofs in a domain-specific language or ``synthetic'' mathematical framework. For instance, if one extends Martin-L\"of's dependent type theory with propositional truncations, function extensionality, uniqueness of identity proofs, and the axiom of choice, one gets the domain-specific language of sets, which is used in \mathlib{}. Similarly, \emph{homotopy type theory} is the domain-specific language of homotopy types or spaces, in which Martin-L\"{o}f's dependent type theory is extended with the univalence axiom and higher inductive types, and uniqueness of identity proofs is dropped. Other domain-specific approaches to formalization include parametrized spectra and synthetic algebraic geometry. We will focus here on homotopy type theory.

The primitive objects in homotopy type theory---that is, the types---are homotopy types or spaces, and everything definable in homotopy type theory is automatically invariant under homotopy. Familiar spaces such as the spheres $S^n$ are defined as higher inductive types, essentially as homotopy pushouts, rather than as subspaces of Euclidean spaces. For each pointed type $A$, we can define the $n$\supth{} homotopy group $\pi_n(A)$, which is an example of a ``homotopy invariant'' construction, as the set-truncation of the space of pointed maps from $S^n$ to $A$. Computations of homotopy groups within the setting of homotopy type theory have to be done entirely through homotopy-invariant methods. 

In particular, there is a construction of an isomorphism $\pi_3(S^2) \cong \ZZ$ in homotopy type theory \cite[Corollary 8.6.19]{hottbook}, which we will describe in more detail in \zcref{sec:math}. This proof mirrors in some ways and departs in others from the ``analytic'' proof in traditional foundations. An advantage of doing homotopy theory ``synthetically'' in this manner is that proofs and constructions apply to all models of homotopy type theory, whereas proofs in point-set topology only apply to point-set topology; in particular, all constructions in homotopy type theory apply to all Grothendieck $\infty$-topoi \cite{shulman-infinity-topos}.

\subsection{Open questions}

One motivation for domain-specific ``synthetic'' formal systems is their amenability to computer formalization \cite{Voevodsky}, and in fact the construction in homotopy type theory of the isomorphism $\pi_3(S^2) \cong \ZZ$ has been formalized multiple times in multiple proof assistants.\footnote{Among the currently actively libraries, see the \href{https://github.com/agda/cubical/blob/master/Cubical/Homotopy/Group/Pi3S2.agda}{Cubical Agda formalization} and the \href{https://github.com/HoTT/Coq-HoTT/blob/master/theories/Homotopy/PiSpheres.v}{Coq-HoTT formalization}. The latter isomorphism was actually added after our experiment was complete, though for the reasons noted below, all of the essential ingredients for this calculation were already present in that library.} However, homotopy type theory libraries tend to be much smaller than Mathlib, providing less training data for an AI. Would an autoformalization agent, with no specific training in homotopy type theory, be able to formalize the homotopy type theory proof of this theorem?

We had a specific test library in mind, which like Mathlib had sufficient development to state the target theorem $\pi_3(S^2) \cong \ZZ$, but was missing the most significant ingredients in its proofs. The \agdaunimath{} library was created in 2021 with the ambition to serve as a general purpose library for formalized mathematics in univalent foundations \cite{agda-unimath}. 

One of our motivations for conducting this experiment is to test autoformalization capabilities in an area of synthetic mathematics. Another motivation is to attempt to address an equally important question of whether autoformalization of mathematics can be done in a responsible way that ultimately supports the libraries that it relies upon. In recent months, there have been occasions where autoformalizers have arguably hurt human-led formalization efforts more than helping them. In an incident described in \cite{avigad2026mathematiciansageai}, a competing autoformalization project piggybacked upon and then stole credit from a human-led collaboration, without contributing much in the way of usable code. Arguably none of the AI-generated code was of ``library quality,'' suitable for being contributed back into \mathlib{}.

\subsection{Overview of results and outline}

The task was for {\codex} to successfully formalize the computation of $\pi_3(S^2)$, with human prompts for instructions but otherwise in an autonomous manner. Importantly, no human was to supply any piece of the code. Secondary, we set the goal of having it produce high-quality code suitable for the {\agdaunimath} library. In particular, our aim was to autoformalize intermediate results in full and not just the portions of those theorems needed for the isomorphism $\pi_3(S^2) \cong \ZZ$.

Initially, we provided {\codex} with a project description, explaining the overall goal of the project and the methods we expected it to be using. Based on this description, we used an autogenerated formalization plan, with intermediate targets. After a little under one month (June 3 to June 25, 2026), {\codex}, guided by the plan and our repeated prompting, was able to complete the formalization plan, producing a formal proof of $\pi_3(S^2) \cong \ZZ$ in homotopy type theory.

All of our code can be found on github, which also contains anonymized versions of our full codex session logs, including our prompts and the model's responses:
\begin{center} 
\href{https://github.com/emilyriehl/Codex-Homotopy-Group/tree/arxiv-v1}{github.com/emilyriehl/Codex-Homotopy-Group/tree/arxiv-v1}
\end{center}

In \zcref{sec:math}, we give a brief introduction to homotopy type theory and then sketch the proof of the calculation $\pi_3(S^2) \cong \ZZ$ in homotopy type theory following \cite{hottbook}. We also give more background about the {\agdaunimath} library.

In \zcref{sec:experiment}, we give more details about our experimental setup, our resource usage, and the results, in particular describing the other formalization milestones achieved along the way. We also highlight bottlenecks and other inefficiencies and describe how they were overcome. 

In \zcref{sec:conclusions}, we record our takeaways from the experiment, list various future technical improvements, and revisit the important discussion of what this autoformalization means for the {\agdaunimath} library.

\subsection{Acknowledgments}

The first, third, and fourth authors are supported by a grant through DARPA's Exponentiating Mathematics program (HR0011262E019). The credits used for this autoformalization experiment were donated by OpenAI through this initiative. 

In addition, the first author is also supported by the NSERC, while the third author is also supported by the NSF (DMS-2507077) and by the AFOSR (FA9550261B013). Help setting up \codex{} was provided by the Johns Hopkins Data Science and AI Institute, which employs the second author. 

While AI was used, with human supervision, to generate a plan for the experiment, generate the experimental logs, and generate the \agda{} code and documentation, this paper describing our results was entirely written by the four human authors listed above.

As described in the text below, our {\codex} autoformalization succeeded because it was able to learn from the human authors of the original homotopy type theory proofs as well as past work by human experts formalizing these results in other homotopy type theory libraries. We have endeavored to provide precise references throughout the text and wish to convey our appreciation for this beautiful mathematics.

\section{Mathematical target and context}\label{sec:math}

\subsection{Basic concepts of homotopy type theory}\label{ssec:hott-intro}

In homotopy type theory, the axiom \emph{uniqueness of identity proofs} is dropped in favor of the \emph{homotopy interpretation of type theory}, where types are spaces\footnote{More precisely, the types in homotopy type theory should be thought of as homotopy types, or $\infty$-groupoids, or anima---or, more generally, as objects in an $\infty$-topos. We adopt the common homotopy theorists' convention of using ``spaces'' for these notions.}, dependent types are fibrations, and the type $a=b$ of identifications between two points $a$ and $b$ of the same type is the type of paths between them \cite{AwodeyWarren2009}. The type $a = b$ exists whenever $a$ and $b$ are points belonging to the same type, but it is not necessarily the case that $a = b$ itself contains any points. We write $p : a = b$ to indicate that $p$ is a path from $a$ to $b$. The path $p$ is a witness that $a$ and $b$ are ``equal,'' identifying $a$ with $b$ in an explicit structural sense, as is consistent with the common mathematical practice of treating two objects as the same if they are the same in a structural sense. In particular, for any $a$, we have $\refl : a = a$, and moreover the family of identity types associated to a type $A$ is freely generated by these ``constant'' paths: it is an inductive type with $\refl$ as its only constructor.

Given two identifications $p,q:a=b$, we might ask whether they are themselves the same. The type-theoretic way of stating that $p$ and $q$ are the same is via the identity type $p=q$, and we see here at once that this process can be iterated, giving rise to a tower of iterated identity types. In particular, we can define the \emph{iterated loop spaces} of a pointed type $(A,a)$ by
\begin{equation*}
  \Omega^0(A,a)\coloneqq(A,a)\qquad\text{and}\qquad \Omega^{n+1}(A,a)\coloneqq\Omega^n(a=a,\refl).
\end{equation*}

When encountering a type in the context of homotopy type theory, one typically wants to characterize its identity type at the earliest occasion. The identity type of the universe is by the \emph{univalence axiom} postulated to be equivalent to equivalences; that is, an identification of types is equivalent to an equivalence between them. This axiom implies that two sets are the same if there is a bijection between them, two algebraic structures such as two groups or two rings are the same if they are isomorphic, two subsets are the same if they contain the same elements, two equivalence relations are the same if they relate the same elements, and so on.

Identity types can also be used in the specification of a new kind of inductive type: \emph{higher inductive types}. The first example of a higher inductive type is the circle, which is a type $S^1$ generated by a point $\pt$ and a loop
\begin{equation*}
  \loopcircle : \pt=\pt.
\end{equation*}
The circle is then equipped with an induction principle, which tells us that in order to construct a section of a fibration over $S^1$, we have to construct the value of this section at $\pt$ and show that it is a fixed point (as stated using identity types) under the transport of $\loopcircle$. A special case of this induction principle is the universal property of the circle, which states that a map $S^1\to X$ is equivalently described as a free loop $(x,\omega)$ in $X$, where $x:X$ and $\omega:x=x$. That is, we have an equivalence
\begin{equation*}
  (S^1 \to X)\simeq \sum_{x:X}x=x.
\end{equation*}
A beautiful argument via the encode-decode method shows how the univalence axiom can be used to obtain an equivalence, and indeed a group isomorphism, $\Omega S^1\simeq \mathbb{Z}$ \cite{LicataShulman2013}.

More generally, we can define the \emph{suspension} of a type as a higher inductive type $\Sigma X$ with two points $\north,\south:\Sigma X$ and a map
\begin{equation*}
  \merid : X \to \north=\south.
\end{equation*}
The universal property of the suspension states that the evaluation map
\begin{equation*}
  (\Sigma X \to Y)\to \left(\sum_{y,z:Y}X \to y=z\right),
\end{equation*}
given by $f\mapsto (f(\north),f(\south),\ap_f\circ\merid)$ is an equivalence. Here, the function 
\begin{equation*}
  \ap_f : x = y \to f(x) = f(y),
\end{equation*} 
called the \emph{action on identifications} of $f$, is defined inductively by $\ap_f(\refl)\coloneqq\refl$.

Given pointed types $A$ and $B$, with implicit basepoints, we can reformulate the universal property of the suspension as the familiar adjunction
\begin{equation*}
  (\Sigma A \to_\ast B)\simeq_\ast (A \to_\ast \Omega B).
\end{equation*}
Here, we use $\to_\ast$ to denote the type of pointed maps, and indeed the equivalence of the pointed universal property of suspension is a pointed equivalence.

By iterating the suspension, we can define the $n$-sphere for each $n$:
\begin{equation*}
  S^{n+1}\coloneqq\Sigma^n S^1.
\end{equation*}
The universal property of the $n$-sphere $S^n$ provides an equivalence
\begin{equation*}
  (S^n\to_\ast B)\simeq_\ast \Omega^n B
\end{equation*}
for any pointed type $B$.

We note that a type $X$ is considered to be a \emph{set}, or homotopically discrete, if every map $S^1\to X$ is constant, in the sense that the diagonal map
\begin{equation*}
  X \to (S^1\to X)
\end{equation*}
given by $x\mapsto (t\mapsto x)$ is an equivalence. In particular, in homotopy type theory we classify only some types as sets. This way, we view homotopy type theory as an \emph{extension} of traditional set-level mathematics.

Using the idea of higher inductive types, we can also define the \emph{set truncation} of a type $X$. The set truncation of $X$ is a type $\|X\|_0$, which is a set according to the previous definition, that comes equipped with a map $X\to\|X\|_0$ that is universal among maps from $X$ to a set. There are many equivalent ways of constructing set truncations, and indeed it is customary to work only with its \emph{specification} once it has been implemented. A fruitful way of thinking about the type $\|X\|_0$ is that it is the type of connected components of $X$.

More generally, we say that a type $X$ is $k$-truncated if every map $S^{k+1}\to X$ is constant in the sense that the diagonal map
\begin{equation*}
  X\to (S^{k+1}\to X)
\end{equation*}
is an equivalence. Taking $S^0\coloneqq\{\pm 1\}$ and $S^{-1}\coloneqq\varnothing$, this definition makes sense for $k\geq -2$. A $(-2)$-truncated type is called \emph{contractible}, and a $(-1)$-truncated type is called a \emph{mere proposition} \cite[Definition 3.3.1]{hottbook} or a \emph{proposition} \cite[Definition 12.1.1]{introhott}. The $k$-\emph{truncations} $\|X\|_k$ can be defined for all $k \geq -2$, as the universal $k$-truncated type equipped with a map $X\to\|X\|_k$. A type $X$ is $k$-\emph{connected} if its $k$-truncation $\|X\|_k$ is contractible, while a map $f : X \to Y$ is $k$-\emph{connected} if its \emph{fibers} $\sum_{x : X} f(x) = y$ over all points $y : Y$ are $k$-connected. For example, the $(k+1)$-sphere is $k$-connected.

The homotopy groups of pointed types $A$ are defined in homotopy type theory as the set truncations of the iterated loop spaces:
\begin{equation*}
  \pi_n(A) \coloneqq \|\Omega^n(A)\|_0.
\end{equation*}
For $n\geq 1$, the group structure on these sets is inherited from the groupoidal structure of the identity type, and the Eckmann--Hilton argument, showing that the homotopy groups are abelian for $n\geq 2$, is given in homotopy type theory via path induction \cite[Theorem 2.1.6]{hottbook}. A type $X$ is $k$-connected if and only if its homotopy groups $\pi_n(X,x)$ are trivial for each $n\leq k$ and each $x:X$. For $k \geq 0$, a map $f : X \to Y$ is $k$-connected if and only if the induced map $\|X\|_0 \to \|Y\|_0$ is an equivalence and, for each $x : X$, the induced map on homotopy groups $\pi_n(X,x) \to \pi_n(Y, f(x))$ is an isomorphism for $n \leq k$ and is surjective when $n = k +1$ \cite[Corollary 8.8.5]{hottbook}.

An important difference between homotopy type theory and point-set topology is that the condition that a type $X$ is $k$-truncated is stronger than the condition that $\pi_n(X,x)$ is trivial for all $n>k$ and all $x:X$. There are models of homotopy type theory in which there is a map that is $k$-connected for all $k$, while not being an equivalence \cites[Counterexample 2.10.1 and Remark 2.10.2]{lurie2003}[\S 8.8]{hottbook}{shulman-infinity-topos}.

\subsection{The calculation of \texorpdfstring{$\pi_3(S^2)$}{π₃(S²)} in homotopy type theory}\label{ssec:hott-calculation}

The calculation of the group $\pi_3(S^2)$ in homotopy type theory has three major ingredients that we explain below:
\begin{enumerate}
\item The construction of the long exact sequence of homotopy groups of a fiber sequence.
\item The construction of a particular fiber sequence $S^1 \to_\ast S^3 \to_\ast S^2$, called the Hopf fibration, which together with the long exact sequence gives isomorphisms $\pi_2(S^2)\cong \pi_1(S^1)$ and $\pi_3(S^3)\cong \pi_3(S^2)$.
\item The Blakers--Massey theorem, which implies the Freudenthal suspension theorem, which gives an isomorphism $\pi_2(S^2)\cong\pi_3(S^3)$. 
\end{enumerate}

A \emph{fiber sequence} $F \to_\ast E\to_\ast B$ in homotopy type theory consists of pointed maps $i : F\to_\ast E$ and $p : E \to_\ast B$ and a pointed homotopy $H$, which is often left implicit, witnessing that the following square commutes and is a homotopy pullback:
\begin{center}
  \begin{tikzcd}[
    alt={
      A commutative square with objects F, E, 1, and B, listed in reading order.
      There are maps from F to E, from E to B, from F to 1, and from 1 to B.
      The map from F to E is the inclusion of the fiber and is labeled i;
      the map from E to B is the fibration and is labeled p.
      The other maps are uniquely determined as pointed maps, and are thus unlabeled.
    }
    ]
    F \arrow[r, "i"] \arrow[d] & E \arrow[d, "p"] \\
    \unit \arrow[r] & B.
  \end{tikzcd}
\end{center}
Every fiber sequence is uniquely determined by a pointed map $p$ of type $E \to_\ast B$, as the fiber, defined above, can be shown to define the homotopy pullback. However, the type of such data can be shown to be equivalent to the type that is given by a pair of pointed maps, together with an equivalence between $F$ and the fiber of $p$ commuting with the inclusions into $E$, and this expanded form of the data is more convenient in situations, such as ours, where we wish to explicitly specify $F$, $E$, and $B$.  The \agdaunimath{} library is very deliberate about representation independent design choices such as these.

Any fiber sequence $F \to_\ast E\to_\ast B$ induces a fiber sequence $\Omega B \to_\ast F \to_\ast E$. Iterating this process, we obtain a tower of types $\cdots\to_\ast \Omega B\to_\ast F \to_\ast E\to_\ast B$. This construction was first formalized by Voevodsky in \cite{Voevodsky2011UnivalentBasics}. Applying the set truncation to this tower yields the \emph{long exact sequence} of homotopy groups~\cites[Theorem 8.4.6]{hottbook}[Theorem 4.3.4]{AvigadKapulkinLumsdaine2015}:

\begin{center}
  \begin{tikzcd}[
    column sep=2em,
    alt={
      A segment of the long exact sequence of homotopy groups associated to
      the fiber sequence F to E to B.
      This segment displays the groups pi sub n plus 1 of B,
      pi sub n of F, E, and B,
      and finally pi sub n minus 1 of F.
      The maps are the connecting homomorphisms together with the maps induced
      by i and p on homotopy groups.
    }
    ]
    \cdots \arrow[r]
    & \pi_{n+1}(B) \arrow[r, "\partial"]
    & \pi_n(F) \arrow[r, "i_\ast"]
    & \pi_n(E) \arrow[r, "p_\ast"]
    & \pi_n(B) \arrow[r, "\partial"]
    & \pi_{n-1}(F) \arrow[r]
    & \cdots.
  \end{tikzcd}
\end{center}

The Hopf fibration $S^1 \to_\ast S^3 \to_\ast S^2$ is a fiber sequence constructed as a special case of the \emph{Hopf construction} associated to a connected $H$-space $A$. An \emph{H-space} is a type $A$ equipped with a binary operation $\mu:A\times A \to A$ and a unit $e:A$ equipped with identifications $\lambda_x:\mu(e,x)=x$ and $\rho_x:\mu(x,e)=x$ for each $x:A$ \cite[Definition 8.5.4]{hottbook}. In the context of homotopy type theory, it is often beneficial to assume an extra coherence condition: an identification $\lambda_e=\rho_e$ \cite{BuchholtzRijke2018}. Such H-spaces are called \emph{coherent} \cite{BuchholtzChristensenFlatenRijke2025}. The associated fiber sequence  $A \to_\ast A \join A \to_\ast \Sigma A$ has fiber $A$, base space the suspension of $A$, and total space $A \join A$ \cite[Lemma 8.5.7]{hottbook}, where the \emph{join} of types $A$ and $B$ is defined to be the homotopy pushout:

\begin{center}
\begin{tikzcd}[
  alt={The defining homotopy pushout square of the join of two types A and B has the type A times B at the top-left corner, a projection into the type A at the lower left corner, a projection into the type B at the top right corner, and the type A join B at the lower left corner.}]
A \times B \ar[d, "\pi_A"'] \ar[r, "\pi_B"] & B \ar[d] \\ A \ar[r] & A \join B. 
\end{tikzcd}
\end{center}

The classical Hopf fibration $S^1 \to_\ast S^3 \to_\ast S^2$ arises when this construction is specialized to the $H$-space $S^1$. Here, the $H$-space structure is given by a multiplication map $S^1 \times S^1 \to S^1$ corresponding to complex multiplication on the unit circle. Such a map can indeed be defined synthetically via the universal property and induction principle of the circle \cite[Lemma 8.5.8]{hottbook}. By the Hopf construction, we obtain a fiber sequence $S^1 \to_\ast S^1 \join S^1 \to_\ast S^2$. Using the equivalence $S^1\ast S^1\simeq S^3$, which requires a proof, we obtain the Hopf fibration. Having constructed the Hopf fibration on the $2$-sphere, we obtain a long exact sequence:

\begin{center}
  \begin{tikzcd}[
    column sep=1.8em,
    alt={
      A segment of the long exact sequence of homotopy groups associated to
      the Hopf fibration.
      This segment displays the groups pi sub n plus 1 of the 2 sphere,
      pi sub n of the 1 sphere, the 3 sphere, and the 2 sphere,
      and finally pi sub n minus 1 of the 1 sphere.
      The maps are the connecting homomorphisms together with the maps induced
      by i and p on homotopy groups.
    }
    ]
    \cdots \arrow[r]
    & \pi_{n+1}(S^2) \arrow[r, "\partial"]
    & \pi_n(S^1) \arrow[r, "i_\ast"]
    & \pi_n(S^3) \arrow[r, "p_\ast"]
    & \pi_n(S^2) \arrow[r, "\partial"]
    & \pi_{n-1}(S^1) \arrow[r]
    & \cdots.
  \end{tikzcd}
\end{center}

Since the circle is $1$-truncated, we have $\pi_n(S^1)\cong 0$ for $n>1$. Furthermore, since $S^2$ is $1$-connected and $S^3$ is $2$-connected, we have that $\pi_1(S^2)\cong 0$ and $\pi_2(S^3)\cong 0$. Thus, we obtain the following exact sequence of groups:
\begin{center}
  \begin{tikzcd}[
    column sep=1.8em,
    alt={
      A segment of the long exact sequence of homotopy groups associated to
      the Hopf fibration.
      This segment displays the groups pi sub n plus 1 of the 2 sphere,
      pi sub n of the 1 sphere, the 3 sphere, and the 2 sphere,
      and finally pi sub n minus 1 of the 1 sphere.
      The maps are the connecting homomorphisms together with the maps induced
      by i and p on homotopy groups.
    }
    ]
      0          \arrow[r]
    & \pi_3(S^3) \arrow[r]
    & \pi_3(S^2) \arrow[r]
    & 0          \arrow[r]
    & 0 \arrow[r]
    & \pi_2(S^2) \arrow[r]
    & \pi_1(S^1) \arrow[r]
    & 0,
  \end{tikzcd}
\end{center}
giving isomorphisms $\pi_3(S^3)\cong \pi_3(S^2)$ and $\pi_2(S^2)\cong\pi_1(S^1)$. In particular, since $\pi_1(S^1)\cong\ZZ$ was proven in homotopy type theory in \cite{LicataShulman2013}, we find that $\pi_2(S^2)\cong\ZZ$.

The computation of $\pi_3(S^2)$ therefore requires one more ingredient: a proof that $\pi_2(S^2)\cong\pi_3(S^3)$. This isomorphism is obtained from the Freudenthal suspension theorem, which asserts that for any pointed $n$-connected type $A$ for $n \geq -1$, the map
\begin{equation*}
  A \to \Omega \Sigma A
\end{equation*}
is $2n$-connected. LeFanu Lumsdaine's beautiful proof was reported in \cite[Theorem 8.6.4]{hottbook} and a formalization was reported in \cite{Licata2013HoTTProgress}. Since the $2$-sphere is $1$-connected, we see that the map
\begin{equation*}
  S^2\to \Omega S^3
\end{equation*}
is $2$-connected, so that it induces an isomorphism $\pi_2(S^2)\cong \pi_3(S^3)$. Thus, the computation $\pi_3(S^2)\cong \ZZ$ in homotopy type theory is completed.

The Freudenthal suspension theorem in homotopy type theory is typically derived as a corollary of the Blakers--Massey theorem. The Blakers--Massey theorem asserts that for a homotopy pushout square
\begin{center}
  \begin{tikzcd}[
    alt={
      A commuting square with objects A, B, C, and D, and named arrows f from A to B and g from A to C, and unnamed arrows from B to D and from C to D.}
    ]
    A \ar[r,"f"] \ar[d,swap,"g"] & B \ar[d] \\
    C \ar[r] & D,
  \end{tikzcd}
\end{center}
if the map $f:A\to B$ is $m$-connected and the map $g:A\to C$ is $n$-connected, then the \emph{gap map}, which is the universal map into the homotopy pullback $A\to B\times_D C$, is $(m+n)$-connected. To obtain the Freudenthal suspension theorem from Blakers--Massey, take $B\coloneqq\unit$ and $C\coloneqq\unit$, so that the homotopy pushout $D$ is the suspension $\Sigma A$ and the homotopy pullback is $\Omega\Sigma A$. The Blakers--Massey theorem was first proven and formalized in homotopy type theory in \cite{HFLL}.

\subsection{The \agdaunimath{} library}

\begin{figure}[!t]
  \centering
  \includegraphics[alt={Dependency graph of the agda-unimath library,
  with nodes colored according to their namespace.},width=\textwidth]{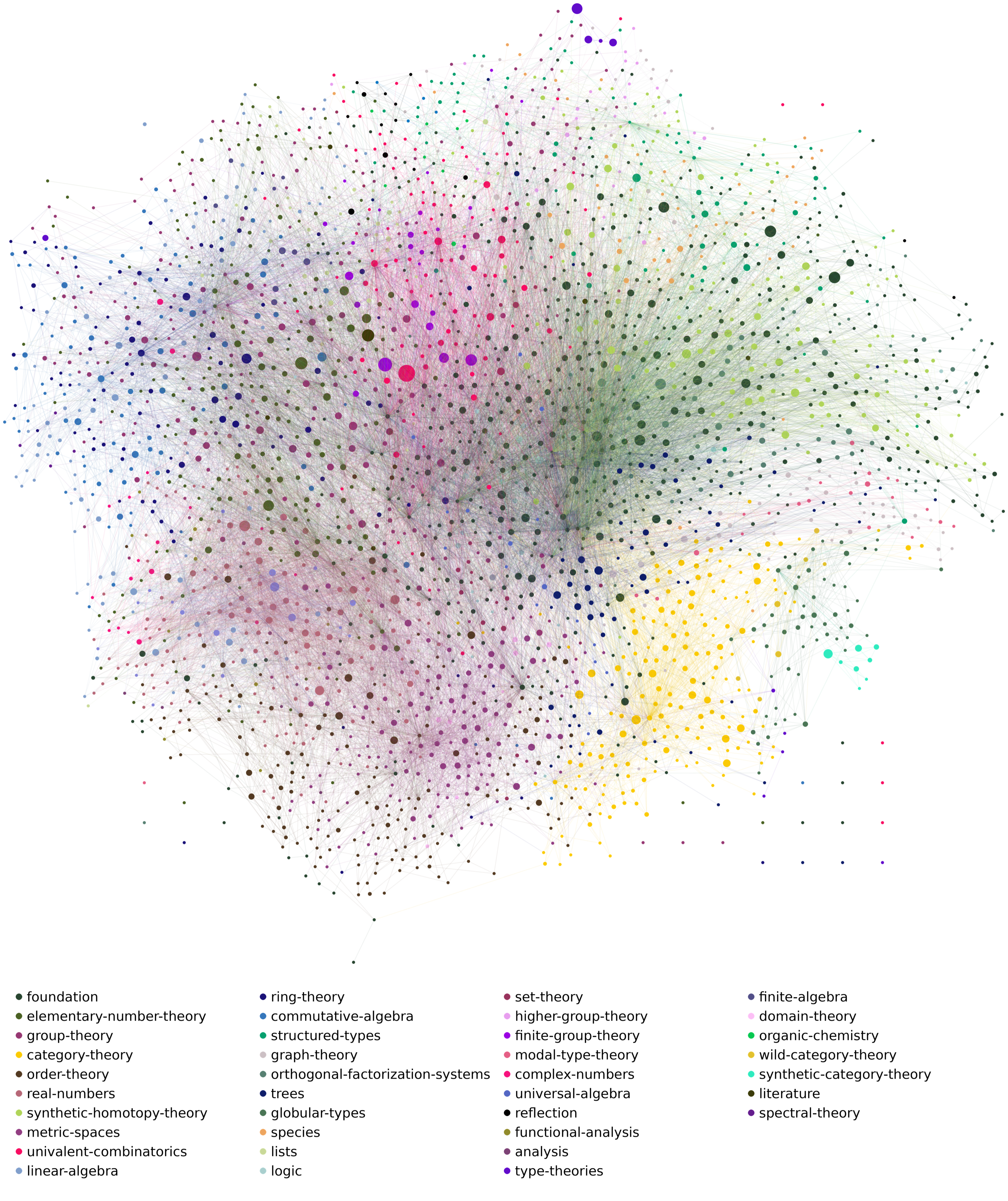}
  \caption{Dependency graph of the \agdaunimath{} library, created by Fredrik Bakke.}
\end{figure}

The \agdaunimath{} library is one of the largest formalization projects for the \agda{} proof assistant. It currently contains over $3000$ pages, with concepts and theorems from a wide variety of mathematical subjects. It is built on the univalent foundations for mathematics, and treats all its subjects from a univalent point of view. Although its foundations are homotopy type theory, synthetic homotopy theory is not the main focus of the library, which made it suitable for our experiment.

The \agdaunimath{} library is organized by mathematical subject. The foundational topics of univalent mathematics are treated in the folders \agdainline|foundation| and \agdainline|foundation-core|. These are organized in two folders for bootstrapping purposes, so as to avoid cyclic dependencies of the modules. The central folders for our experiment are: 
{
  \AgdaCodePlain
  \begin{agdacode}
    foundation
    foundation-core
    group-theory
    structured-types
    synthetic-homotopy-theory
  \end{agdacode}
}
Within the foundation folder, one finds modules such as \agdainline|identity-types|, \agdainline|equivalences|, and \agdainline|truncated-types|. The library is designed around a \emph{one-concept-per-file} paradigm. The organization of each file is modeled after $n$Lab pages: Modules start with an \emph{Idea}-section, explaining in natural languare the idea of a concept or named entity, then proceeds with one or several definitions of that concept, and afterwards proves basic properties about it. References and pointers to related files may be given at the end of the module. The \agdaunimath{} libary uses sharp, disambiguative conventions for concept declarations, where everything that can be named or disambiguated in natural language gets its own file. As an example, we show an excerpt from the module \agdainline|synthetic-homotopy-theory.circle|, in which the standard circle is postulated.

\begingroup
\AgdaCodeName{plain}{with}
\AgdaCodeName{plain}{synthetic-homotopy-theory.loop-spaces}
\begin{agdacode}
## Idea

**The circle** is the initial type equipped with a base point and a
[loop](synthetic-homotopy-theory.loop-spaces.md).

## Postulates

```agda
postulate
  𝕊¹ : UU lzero

postulate
  base-𝕊¹ : 𝕊¹

postulate
  loop-𝕊¹ : base-𝕊¹ ＝ base-𝕊¹

free-loop-𝕊¹ : free-loop 𝕊¹
free-loop-𝕊¹ = base-𝕊¹ , loop-𝕊¹

𝕊¹-Pointed-Type : Pointed-Type lzero
𝕊¹-Pointed-Type = 𝕊¹ , base-𝕊¹

postulate
  ind-𝕊¹ : induction-principle-circle free-loop-𝕊¹
```
\end{agdacode}
\endgroup

Here, the entry \agdainline|induction-principle-circle| is defined in the module 
\begin{agdacode}
synthetic-homotopy-theory.universal-property-circle
\end{agdacode} 
which specifies what \emph{a} circle is in homotopy type theory without declaring the standard circle. The infrastructure written for circles is written representation independently so as to be applicable to any type that is a circle, for instance, the type $S^0\ast S^0$.

Finally, we note that \agdaunimath{} formalizes the theory of groups in two parallel ways. The ordinary notion of group is the usual one, where a group is a set with a binary operation and a unit element, satisfying the axioms of a group. We also formalize the notion of \emph{concrete group}, which is a pointed connected $1$-type. The idea is that identifications of an element with itself can be thought of as the \emph{symmetries} of that element. For example, the type of $n$-element types is the concrete group presenting the $n$\supth{} symmetric group $S_n$, since the identifications of the standard $n$-element type $\Fin{n}$ are bijections $\Fin{n}\cong \Fin{n}$. Thus, a concrete group \emph{concretely} presents a group as the group of symmetries of an object. Given a pointed connected $1$-type $BG$, the underlying type of the group is $G:=\Omega BG$. From the perspective of concrete groups, a group homomorphism $G\to H$ is a pointed map $BG\to_\ast BH$. All of group theory can be developed in this way \cite{symmetry}. The idea of treating group theory via classifying spaces in the context of homotopy type theory was first explored in \cite{shulman2015} and a first systematic study was reported in \cite{BvDR}.

At the moment of our experiment, \agdaunimath{} defined the \emph{concrete} homotopy groups as follows:

\begin{agdacode}
### The concrete homotopy groups

```agda
module _
  {l : Level} (n : ℕ) (A : Pointed-Type l)
  where

  concrete-homotopy-group : Concrete-Group l
  concrete-homotopy-group =
    concrete-group-Pointed-Type (iterated-loop-space n A)
```
\end{agdacode}

Here, the concrete group associated to a pointed type $(A,a)$ is the $1$-truncation of the connected component of the base point $a$. Since the underlying set of a concrete group $BG$ is its loop space, we see here that \agdainline|concrete-homotopy-group n A| is the classifying space of $\pi_{n+1}(A)$. This unfortunate shift in indexing will come to play a part in our experiment.

\section{Design and results of the experiment}\label{sec:experiment}

In this section, we describe the design and results of the autoformalization experiment, which was initiated on June 3\suprd{}, 2026; the main theorem was completed on June 25\supth{}, 2026.

\subsection{Experimental design}

An ideal test of autoformalization capabilities in homotopy type theory involves a theorem that is stateable using the previously-defined infrastructure of an ambient library, but which requires a significant amount of work to prove. The calculation of $\pi_3(S^2) \cong \ZZ$ perfectly fit this desiderata. The \agdaunimath{} library already contained definitions of the spheres and homotopy groups, the notion of a group isomorphism, and the group of integers. But it was missing most of the essential ingredients in this calcuation: in particular the long exact sequence of a fibration, the Hopf fibration, and the Freudenthal suspension theorem.

After we identified our formalization target, we checked to see whether any human-led teams were working on any of these components. On a Discord server associated with the \agdaunimath{} library, we wrote 
\begin{quote}
  One of the ways we're trying to be responsible in experimenting with AI for formalization is to check first whether any of our autoformalization targets are currently works in progress by humans who should be given priority. To that end is anyone currently working on formalizing the Hopf fibration or the calculation of $\pi_3(S^2)$ for \agdaunimath{}? This also involves the long exact sequence of the fibration. \ldots right now we want to make sure that the mere existence of this code does not disrupt/conflict with/undermine any human efforts.
\end{quote}
After we determined that there was no known formalization effort in this direction, we settled on this project.

As none of the first, third, and fourth authors had experience using a coding agent, we submitted a request for research software engineering support from the Johns Hopkins Data Science and AI Institute  (DSAI), which put them in contact with the second author. 

Using an email that the third author wrote to describe the project to colleagues at the DSAI as a prompt (preserved in the repository as email.md\footnote{\href{https://github.com/emilyriehl/Codex-Homotopy-Group/blob/arxiv-v1/email.md}{github.com/emilyriehl/Codex-Homotopy-Group/blob/arxiv-v1/email.md}}), {\claude} Code generated a guiding formalization plan document\footnote{\href{https://github.com/emilyriehl/Codex-Homotopy-Group/blob/arxiv-v1/FORMALIZATION-PLAN.md}{github.com/emilyriehl/Codex-Homotopy-Group/blob/arxiv-v1/FORMALIZATION-PLAN.md}} that correctly identified a subtlety in the \agdaunimath{} statement of the main goal:
\begin{quote}
Important indexing note: \agdainline|synthetic-homotopy-theory.homotopy-groups| documents
an ``Obi-wan error'' in \agdainline|concrete-homotopy-group|: the index \agdainline|n| names the
\agdainline|(n+1)|\supst{} abstract homotopy group. Therefore $\pi_3(S^2)$ should correspond to
\agdainline|concrete-homotopy-group 2 (sphere-Pointed-Type 2)|, not index \agdainline|3|. Candidate final shape:

\begin{agdacode}
iso-Group
  ( group-Concrete-Group
    ( concrete-homotopy-group 2 (sphere-Pointed-Type 2)))
  ( ℤ-Group)
\end{agdacode}
\end{quote}

A similar formalization plan produced by \codex{} failed to notice the off-by-one error, stating the theorem incorrectly as an isomorphism 
\begin{agdacode}
iso-Concrete-Group (concrete-homotopy-group 3 sphere-2) integer-Concrete-Group 
\end{agdacode}
Thus, we initiated \codex{} on the {\claude} Code-generated formalization plan.

We set up \codex{} in a shared private github repository so that the prompting could be done by any of the first, third, or fourth authors using multiple computers.

Before attempting any autoformalization, we provided \codex{} with skills and reference documentation, instructing it in the importance of using \agda{} \agdainline|--without-K|, which is to say in homotopy type theory mode, and matching the \agdaunimath{} \href{https://unimath.github.io/agda-unimath/DESIGN-PRINCIPLES.html}{library design principles} and \href{https://unimath.github.io/agda-unimath/CODINGSTYLE.html}{style guides}.

The skills are collections of instruction and reference documents that the agent is told to consult as it works, whose initial versions were generated using {\claude} Code from the \agdaunimath{} source and documentation. The reference skill contains four documents describing:
\begin{itemize}
\item the writing conventions of the library;
\item how to search the library's more than 3000 files for an existing lemma;
\item the foundational vocabulary of the library: universes, $\Sigma$-types, identity types, equivalences, and truncations;
\item and how to interpret the \agda{} error messages that arise when working in this mode.
\end{itemize}

The workflow skill tells the agent how to work: search the library before proving a lemma from scratch, make small edits and typecheck often, never introduce a postulate, and keep a status report and a chat log up to date as a record of the experiment. During the experiment, \codex{} extended this material at our prompting, for instance adding a fifth reference document on formalization practice in homotopy type theory\footnote{\href{https://github.com/emilyriehl/Codex-Homotopy-Group/blob/arxiv-v1/.codex/skills/agda-unimath-reference/references/hott-skills.md}{github.com/emilyriehl/Codex-Homotopy-Group/blob/arxiv-v1/.codex/skills/agda-unimath-reference/references/hott-skills.md}}.

We also had to decide in advance how formalized code would be judged. Testing before the project began had shown that the \agda{} MCP server, an MCP server that gives coding agents interactive access to \agda{} (discussed further below), would report ``ok-complete'' for files it had merely scope-checked, passing code that \agda{} itself rejects. The repository therefore includes a short script, check.sh\footnote{\href{https://github.com/emilyriehl/Codex-Homotopy-Group/blob/arxiv-v1/check.sh}{github.com/emilyriehl/Codex-Homotopy-Group/blob/arxiv-v1/check.sh}}, that runs the actual \agda{} binary on a file with the library flags applied and unsolved metavariables disallowed. The script requires no network access, so \codex{} could run it from inside its own sandbox. A module counted as proven only when check.sh ran cleanly and no holes remained in the file; the instructions to the agent state that nothing else would be considered success.

\subsection{Initial autoformalization}

The initial formalization plan noted 35 dependencies, including the main theorem, 11 of which were missing: 
\begin{enumerate}
  \item the definition and basic theory of pointed fiber sequences,
  \item the long exact sequence of homotopy groups from a pointed fiber sequence,
  \item the isomorphism that arises from a four term exact sequence with zeros at its endpoints,
  \item the Hopf fibration $A \to_\ast A \ast A \to_\ast \Sigma A$ associated to a connected $H$-space $A$,
  \item the connected $H$-space structure on $S^1$,
  \item the Hopf fibration $S^1 \to_\ast S^3 \to_\ast S^2$ arising from equivalences $S^3 \simeq S^1 \ast S^1$ and $S^2 \simeq \Sigma S^1$,
  \item the computation $\pi_n S^1 \cong 0$ for $n > 1$,
  \item the Freudenthal suspension theorem as a corollary of the Blakers--Massey theorem (which was also missing),
  \item the isomorphisms $\pi_{n+3} S^{n+3} \cong \pi_{n+2} S^{n+2}$ for $n \geq 0$ that follow from Freudenthal, together with a separate base case\footnote{The formalization plan erroneously lists $\pi_3 S^3 \cong \pi_2 S^2$ as the required base case, but this is in the range of Freudenthal. The actually required base case is instead $\pi_2S^2 \cong \pi_1 S^1$, which follows from the long exact sequence of the Hopf fibration as noted above. {\codex} autoformalized this by showing that $S^3 \simeq \Sigma \Sigma S^1$ is 2-connected, so its homotopy groups in degrees $1$ and $2$ vanish.} $\pi_2 S^2 \cong \pi_1 S^1$,
  \item the computation $\pi_n S^n \cong \ZZ$ for $n > 1$ (with the $n=1$ case already in the library),
  \item the deduced isomorphism $\pi_3(S^2) \cong \ZZ$ arising from the exact sequence \[ 0 \cong \pi_3(S^1) \to \pi_3(S^3) \to \pi_3(S^2) \to \pi_2 S^1 \cong 0. \]
\end{enumerate} 

The initial prompting on June 3\suprd{} targeting the first of these --- the definition of general pointed fiber sequences --- was done as a group with the code carefully reviewed. We requested a repackaging of the \codex{}-generated \agdainline|fiber-sequence-Pointed-Type| including a renaming of some of its fields. We also directed \codex{} to use some existing API for the coherence involved in a commuting triangle of pointed maps.

The next objective, to formalize the long exact sequence of homotopy groups associated to a fiber sequence of pointed types, was the first major test.

\subsection{Progress reports}

Over the subsequent weeks, we made a number of errors associated to new users:
\begin{itemize}
\item After a few days, we learned how to change the model reasoning level to gpt-5.5 xhigh, which was used for the remainder of the experiment.
\item After a few more days, we learned how to set up an API key associated with the organization which had the donated credits. This is why our project spending report dates from June 9\supth{} rather than the project start date of June 3\suprd{}. 
\item The following day, we realized that {\codex} was not configured to use the \href{https://libraries.io/npm/agda-mcp-server}{agda-mcp-server} we had identified. Over the next few days, we figured out that this needed to be configured on each individual machine interacting with the git repo.
\end{itemize}

These last two improvements were made while the third author attending a workshop on ``Formalization of Higher Categories'' at the Mittag-Leffler institute and got some tips from colleagues. Andrej Bauer explained the use of plan mode to generate more complicated prompts and also suggested instructing \codex{} to attempt some top-down reasoning, decomposing the target isomorphism in terms of its easily-stable prerequisites, with holes left in the places where the formalized proofs could not be easily defined.

It is hard to say how much the MCP server helped, as we made no controlled comparison, and it was configured in the same week as the other corrections just described. We do know how much it was used, and how often it misled us. At the end of the experiment we asked \codex{} to extract every \agda{} MCP call from its own session logs. There were 711 of them, so the server was used constantly, mostly to load files, inspect goals, and infer types. But the extraction also flagged 290 instances of misleading behavior, including 70 replies that reported success while containing an \agda{} error, and a handful in which an editing tool claimed a goal was solved although the edited file no longer typechecked. The same warning recurs throughout the session logs: check.sh, and not the server, has the final word, because the server could report success on code that \agda{} rejects. The report we prepared for the server maintainers\footnote{\href{https://github.com/emilyriehl/Codex-Homotopy-Group/tree/arxiv-v1/agda-mcp-ux-report}{github.com/emilyriehl/Codex-Homotopy-Group/tree/arxiv-v1/agda-mcp-ux-report}} is discussed in \zcref{ssec:future-improvements}.

\subsection{Bottlenecks}\label{ssec:bottlenecks}

A major source of inefficiency was caused by the fact that the primary prompter, the third author, never learned how to reliably convince \codex{} to keep working without her direct supervision. Thus, very little was accomplished on days during which she was occupied by other tasks.\footnote{After the completion of this experiment, the authors learned about goal mode.}

There were various stages where the formalization targets seemed too hard for \codex{} without further expert guidance concerning how to handle the ``invisible mathematics'' that arises in computer formalization, which must concern itself with a level of detail that is absent in the natural-language literature. The most effective tool to overcome such instances involved pointing \codex{} to existing human formalizations found in the Coq-HoTT library \cite{coqhott} using \rocq{}. Once \codex{} learned to use this library as a guide, it autoformalized a third hard theorem, the Blakers--Massey theorem, without our help.

The most significant difficulty arose at various points in the proof of exactness of the long exact sequence of homotopy groups defined from a fiber sequence of pointed types $F \to_\ast E \to_\ast B$. Without further guidance, \codex{} tended to pursue direct proofs of the next targets via messy path algebra that would ultimately fail, and then return unrelated incremental progress, leading the third author to prompt on July 16\supth{}:
\begin{quote}
We seem to get stuck at the same point. Can some of the other homotopy type theory repositories that have formalized the long exact sequence help provide inspiration?
\end{quote}

Despite this generic guidance, \codex{} had difficulties formalizing the exactness of the higher stages of the long exact sequence, attempting to use multiple constructions of the map $\Omega^{n+1} B \to_\ast \Omega^n F$ and struggling to prove that they were homotopic, occasionally due to incompatibilities in sign conventions used either to define the iterated boundary maps or to identify the fibers of the iterated fiber sequences. 

Success was finally achieved by more firmly instructing \codex{} to follow the approach of the Coq-HoTT library: from a general fiber sequence $F \to_\ast E\to_\ast B$, one can construct shifted fiber sequences $\Omega B \to_\ast F \to_\ast E$ and $\Omega E \to_\ast \Omega B \to_\ast F$. We instructed \codex{} to stay the course with a prompt on June 19\supth{}:
\begin{quote}
I'd prioritize the harder upfront equivalence. The goal is not just to have a proof but the right proof that we will ultimately upstream to the agda unimath library. Please update your skills and general references to reflect this ultimate priority.
\end{quote}
It took several days of repeated prompting to ``work very hard'' before this was completed on June 21\supst{}, at which point enough of the long exact sequence was in place to fit into the scaffold (though the results were not yet ``library quality'').

On June 21\supst{}, \codex{} began work on the Hopf fibration. The \agdaunimath{} library already had definitions of an $H$-space and a construction of the coherent multiplication operation on the circle, so \codex{} began by packaging this into an $H$-space structure on the circle and then defining the fiber sequence $A \to_\ast A \ast A \to_\ast \Sigma A$ associated to a general $H$-space $A$. As the spheres are defined in \agdaunimath{} by iterated suspension, identifying the sequence $S^1 \to_\ast S^1 \ast S^1 \to_\ast \Sigma S^1$ with a fiber sequence $S^1 \to_\ast S^3 \to_\ast S^2$ only required an equivalence $S^1 \ast S^1 \simeq S^3$, which proved to be surprisingly difficult.

One way to understand this equivalence is to redefine the spheres as iterated join powers of the 2-point type $\Fin{2}$. For this, \codex{} supplied an inductive proof that the $(n+1)$\supst{} join power of $\Fin{2}$ is equivalent to $S^n$ by demonstrating a more general equivalence $\Sigma X \simeq X \ast \Fin{2}$. Thus, the desired equivalence could be understood as an instance of ``join associativity'', an equivalence between $((\Fin{2}\ast\Fin{2})\ast \Fin{2})\ast\Fin{2}$ and $(\Fin{2}\ast\Fin{2})\ast (\Fin{2}\ast\Fin{2})$.

Rather than tackle this directly, \codex{} insisted on pursuing various side quests until directed more clearly to stop. The prompter finally came up with the right combination of words to get \codex{} to persist in its work, though this was to no avail. A first attempt at a direct comparison equivalence $A \ast (B \ast C) \simeq (A \ast B) \ast C$ stalled after several hours. The next day, \codex{} was asked to pivot to arguing by universal properties, though again this was unsuccessful. The following day, after this failed, the third author prompted \codex{} to look at the Coq-HoTT formalization, which compared both $A \ast (B \ast C)$ and $(A \ast B) \ast C$ to a natively-defined triple join. This took two days of work to implement but was completed successfully on June 25\supth{}.

At this point, the only remaining result was the Blakers--Massey theorem and its corollary, the Freudenthal suspension theorem. While \codex{} initially proposed a direct proof of the isomorphism $\pi_3S^3 \cong \pi_2S^2$ needed for our specific goal, the third author redirected it with a prompt:
\begin{quote}
I've changed my mind. Give me a plan instead that goes after a library quality Freudenthal suspension theorem for general use in agda-unimath together with the special case needed here. What is needed for this?
\end{quote}
Over seven hours of repeated prompting, culminating in the instruction ``Keep working hard and don't stop until the Blakers Massey theorem is proven,'' \codex{} successfully autoformalized a proof of the Blakers--Massey theorem, guided primarily by the Coq-HoTT library, in 5000+ lines of code. 

This final step was started and completed on June 25\supth{}, and the unconditional isomorphism $\pi_3S^2 \cong \ZZ$ was autoformalized 30 minutes later, after a handoff to a different machine. Afterwards \codex{} was asked to prove that $\pi_n S^n \cong \ZZ$ for all positive $n$, and work on polishing up the long exact sequence of the fibration, for instance extending the exactness of pointed sets on the path components $\pi_0 F \to \pi_0 E \to \pi_0 B$. This final step of expanding the development of the long exact sequence took place on June 26\supth{} and June 27\supth{}.

\subsection{Metrics}

The formalization comprises 87 literate \agda{} files and 31,464 lines under the project's source directory, of which about 28,000 are \agda{} code and the rest surrounding mathematical prose. Nearly all of it is contained in the synthetic homotopy theory namespace. The two largest files are the join associativity development (7,436 lines), which was the hardest-fought result of \zcref{ssec:bottlenecks}, and the Blakers--Massey theorem (5,474 lines), which was completed surprisingly quickly, most likely due to the learned route of consulting the Coq-HoTT library for guidance and marginal improvements in user prompting. The repository records 151 commits between June 3\suprd{} and June 27\supth{}, when the autoformalized \agda{} code was written. Later commits were made in preparing the repository for public release, generating the \agda{} MCP user experience report, fixing bugs in the check.sh file, cleaning up the documentation concerning the long exact sequence, and adding anonymized versions of the raw codex logs.

{\setlength{\emergencystretch}{2em}
The total amount of tokens used during June 9\supth{} -- June 27\supth{} by the third author was 2,079,055,195, for a cost of \$1,750.49. As noted above, the prompting done from June 3\suprd{} -- June 8\supth{} was inadvertently done using a free chatGPT account, from which precise usage records cannot be recovered. This period included the first 11 commits on June 3\suprd{}, three commits on June 4\supth{}, and seven commits on June 7\supth{} (including two by the fourth author, who similarly failed to set up the intended API key).\par}

\section{Conclusions}\label{sec:conclusions}

We conclude by offering some lessons from our experiment and also describing plans for further work.

\subsection{Takeaways}

We credit the success of this autoformalization experiment to various sources.

Firstly, the \agdaunimath{} library contained a lot of useful background homotopy theory and group theory, including the calculation that $\pi_1S^1 \cong \ZZ$ together with various other infrastructure related to the circle, as well as general results about loop spaces, suspensions, joins, pushouts, descent, truncation levels, and infrastructure relating abstract and concrete groups. In addition, it had clearly written guides for contributors and an easily searchable library with a uniform style.

Secondly, there was a clearly-written natural language source provided by \cite{hottbook} for the material missing from the \agdaunimath{} library, written using then novel terminology that has since become standard for the field. This appeared to support natural language reasoning to provide outlines for the general homotopy type theoretic proofs. These first two ingredients were the ones that led the authors to believe that this result was within scope of \codex{}.

We did not anticipate how important the existing formalizations in the Coq-HoTT library were to overcoming the bottlenecks described above. Given how central these human-authored guides appeared to be for the long exact sequence, the join associativity, and the proof of Blakers--Massey, we now believe that success of this experiment was ultimately predicated on the fact that humans have done this before. Thus, we feel that it is important to credit the human authors of the formalizations that guided \codex{}. 

Mike Shulman formalized the long exact sequence in Coq-HoTT via what \codex{} considered the ``hard structural route'' of packaging the boundary map $\Omega B \to_\ast F$ of a fiber sequence $F \to_\ast E \to_\ast B$ as a map in its own iterated fiber sequences.

Dan Christensen formalized the join associativity in Coq-HoTT using several tricks that \codex{} copied, including formalizing the final equivalence as a composite involving two instances of join symmetry surrounding a ``twisting'' equivalence:
\[ (A \join B) \join C \simeq C \join (A \join B) \simeq A \join (C \join B) \simeq A \join (B \join C).\]

Mike Shulman formalized Blakers--Massey in Coq-HoTT following the proof of \cite{HFLL, ABFL}, but using a novel construction that he describes in the pull request\footnote{\href{https://github.com/HoTT/Coq-HoTT/pull/1078}{github.com/HoTT/Coq-HoTT/pull/1078}}. Again, {\codex} copied this construction, specializing its work from an arbitrary modality to the truncation modalities.

Our initial impression is that \codex{} is capable of translating formalized source code from one library to another. The incompatibility of formalized mathematics across different libraries has been an important practical problem within the community of formalizers. The capacity of \codex{} to transfer code from one library to another seems to be a major step forward in this regard.

\subsection{Future technical improvements}\label{ssec:future-improvements}

The report we prepared for the maintainers of the \agda{} MCP server, included in the project repository together with its extraction script and the redacted evidence for every suspicious call, traces the most serious failures to a single design problem: the server's replies conflate three separate questions --- whether the request was processed, whether \agda{} accepted the command, and whether any goals, metavariables, or constraints remain open. An agent that takes such a reply at face value can assert a proof that \agda{} rejects. As our attempts to contact the maintainer of that repository have failed, we plan to develop an improved server, tested against the defects catalogued in the report, with this project's formalization as its first test corpus. If interactive feedback can be made trustworthy in this way, a full run of \agda{}, as in check.sh, would be needed only as a final safeguard.

Such a server would matter beyond convenience because of the way \agda{} is used. An \agda{} proof is a program, ordinarily written incrementally through typed holes: a partially written proof can be checked with its gaps left open, and at each gap \agda{} reports the goal and the available hypotheses. A proof advances by inspecting a hole, splitting it into cases, or filling it, and the MCP server exposes to an agent exactly the interaction points that a human sees in the editor. A faithful server would therefore give a human and an agent a shared, machine-checked view of the same unfinished proof, in which either one can leave a hole for the other to attempt --- the top-down style of collaboration suggested to us by Andrej Bauer. The division of labor would then be visible in the proof state itself and checked by \agda{} at every step, rather than coordinated through instructions to the agent. We consider this one of the more promising directions to emerge from the experiment.

A second improvement would automate the prompting itself: a driver that hands the agent one module of the formalization plan at a time and accepts each result only after independent verification, so that progress would not depend on the availability of a human prompter (\zcref{ssec:bottlenecks}). We leave this to future work.

A further place we would like to improve is by giving the model better coaching regarding its communications with us. For instance, {\claude} Code also generated a section on ``Tractability and risks'', describing potential pitfalls of the project. 
It was able to identify sections of the formalization that would be high-effort due to a large amount of definitions or framework missing from {\agdaunimath}.
Curiously, it did so using a mode of speaking not traditionally found among mathematical prose. 
To highlight one such example:

\begin{quote}
  ``\textbf{Hopf fibration risk}: Constructing the fibration from the H-space circle and proving the total space is S{\^{}}3 will stress pushout, join, and fiberwise-descent infrastructure.''
\end{quote}

We also found that {\codex} defaulted to use a similar mode of speaking in the logs we asked it to maintain. Since we are not accustomed to such writing, it was much harder for us to understand, and we did not find the various AI-generated logs to be as useful as we would have liked.

\subsection{Library impacts}

As noted in the introduction, we do not consider an autoformalization project to be completed until it makes a commensurable contribution to the supporting library, in this case \agdaunimath{}. The raw \codex{} generated code in our repository does not yet meet the standards for an \agdaunimath{} PR, much less the standards of the library. It is  insufficiently documented, insufficiently explained, and far too verbose. The formalization code also needs thorough review by us first, and later by library maintainers.

One motivation for writing this report is to facilitate communication with the maintainers of the \agdaunimath{} library. We hope to discuss the possibility of extracting library-quality code to contribute to the library over a series of PRs in a way that is sensitive to the amount of human labor that the review process will require, in particular by contributing our own efforts to the refereeing process. As we view this to be of paramount importance to the evolving interactions between generative AI and computer formalization of mathematics, we plan to update this preprint with the results of our planned future conversations with the broader \agdaunimath{} community.

\begingroup
\clubpenalty=10000
\widowpenalty=10000
\printbibliography[heading=bibintoc,title=References]

@article{ABFL,
author = {Anel, Mathieu and Biedermann, Georg and Finster, Eric and Joyal, André},
title = {A generalized Blakers–Massey theorem},
journal = {Journal of Topology},
volume = {13},
number = {4},
pages = {1521-1553},
doi = {10.1112/topo.12163},
url = {https://londmathsoc.onlinelibrary.wiley.com/doi/abs/10.1112/topo.12163},
year = {2020}
}

@misc{avigad2026mathematiciansageai,
  author        = {Avigad, Jeremy},
  title         = {Mathematicians in the age of AI}, 
  year          = {2026},
  eprint        = {2603.03684},
  archivePrefix = {arXiv},
  primaryClass  = {math.HO},
  url           = {https://arxiv.org/abs/2603.03684}
}

@article{AvigadKapulkinLumsdaine2015,
  author       = {Avigad, Jeremy and Kapulkin, Krzysztof and LeFanu Lumsdaine, Peter},
  title        = {Homotopy Limits in Type Theory},
  journal      = {Mathematical Structures in Computer Science},
  volume       = {25},
  number       = {5},
  pages        = {1040--1070},
  year         = {2015},
  doi          = {10.1017/S0960129514000498},
  eprint       = {1304.0680},
  archivePrefix= {arXiv}
}

@article{AwodeyWarren2009,
  author       = {Awodey, Steve and Warren, Michael A.},
  title        = {Homotopy Theoretic Models of Identity Types},
  journaltitle = {Mathematical Proceedings of the Cambridge Philosophical Society},
  date         = {2009},
  volume       = {146},
  number       = {1},
  pages        = {45--55},
  doi          = {10.1017/S0305004108001783},
  eprint       = {0709.0248},
  eprinttype   = {arxiv},
}

@book{Benthem1979Checking,
  author        = {Benthem Jutting, L. S. van},
  title         = {Checking {L}andau's {\it {G}rundlagen}\ in the {AUTOMATH} system},
  series        = {Mathematical Centre Tracts},
  volume        = {83},
  publisher     = {Mathematisch Centrum, Amsterdam},
  year          = {1979},
  pages         = {iv+120},
  isbn          = {90-6196-147-5},
  mrclass       = {68A40 (02D99)},
  mrnumber      = {664239},
  mrreviewer    = {I.\ Kramosil}
}

@online{symmetry,
  author  = {Bezem, Marc and Buchholtz, Ulrik and Cagne, Pierre
             and Dundas, Bjørn Ian and Grayson, Daniel R.},
  title   = {Symmetry},
  date    = {2026-08-13},
  url     = {https://github.com/UniMath/SymmetryBook/tree/10bb0f6},
  note    = {Commit 10bb0f6},
}

@article{BuchholtzChristensenFlatenRijke2025,
  author       = {Buchholtz, Ulrik
                  and Christensen, J. Daniel
                  and Flaten, Jarl G. Taxer{\aa}s
                  and Rijke, Egbert},
  title        = {Central H-Spaces and Banded Types},
  journaltitle = {Journal of Pure and Applied Algebra},
  volume       = {229},
  number       = {6},
  eid          = {107963},
  date         = {2025},
  doi          = {10.1016/j.jpaa.2025.107963},
  eprint       = {2301.02636},
  eprinttype   = {arxiv},
}

@inproceedings{BvDR,
  author    = {Buchholtz, Ulrik and van Doorn, Floris and Rijke, Egbert},
  title     = {Higher Groups in Homotopy Type Theory},
  booktitle = {Proceedings of the 33rd Annual ACM/IEEE Symposium on Logic in Computer Science},
  series    = {LICS '18},
  year      = {2018},
  pages     = {205--214},
  publisher = {Association for Computing Machinery},
  address   = {New York, NY, USA},
  doi       = {10.1145/3209108.3209150},
}

@article{BuchholtzRijke2018,
  author       = {Buchholtz, Ulrik and Rijke, Egbert},
  title        = {The Cayley--Dickson Construction in Homotopy Type Theory},
  journal      = {Higher Structures},
  volume       = {2},
  number       = {1},
  pages        = {30--41},
  year         = {2018},
  eprint       = {1610.01134},
  archivePrefix= {arXiv}
}

@inproceedings{gonthier:feit-thompson,
  TITLE = {{A Machine-Checked Proof of the Odd Order Theorem}},
  AUTHOR = {Gonthier, Georges and Asperti, Andrea and Avigad, Jeremy and Bertot, Yves and Cohen, Cyril and Garillot, Fran{\c c}ois and Le Roux, St{\'e}phane and Mahboubi, Assia and O'Connor, Russell and Ould Biha, Sidi and Pasca, Ioana and Rideau, Laurence and Solovyev, Alexey and Tassi, Enrico and Th{\'e}ry, Laurent},
  URL = {https://inria.hal.science/hal-00816699},
  BOOKTITLE = {{LNCS}},
  ADDRESS = {Rennes, France},
  EDITOR = {Sandrine Blazy and Christine Paulin and David Pichardie},
  PUBLISHER = {{Springer}},
  SERIES = {LNCS},
  VOLUME = {7998},
  PAGES = {163-179},
  YEAR = {2013},
  MONTH = Jul,
  DOI = {10.1007/978-3-642-39634-2_14},
  HAL_ID = {hal-00816699},
  HAL_VERSION = {v1},
}

@inproceedings{coqhott,
author = {Bauer, Andrej and Gross, Jason and Lumsdaine, Peter LeFanu and Shulman, Michael and Sozeau, Matthieu and Spitters, Bas},
title = {The HoTT library: a formalization of homotopy type theory in Coq},
year = {2017},
isbn = {9781450347051},
publisher = {Association for Computing Machinery},
address = {New York, NY, USA},
doi = {10.1145/3018610.3018615},
booktitle = {Proceedings of the 6th ACM SIGPLAN Conference on Certified Programs and Proofs},
pages = {164–172},
numpages = {9},
location = {Paris, France},
series = {CPP 2017}
}

@article{hales:kepler,
  author = {
    Hales, Thomas
    and Adams, Mark
    and Bauer, Gertrud
    and Dang, Tat Dat
    and Harrison, John
    and Hoang, Le Truong
    and Kaliszyk, Cezary
    and Magron, Victor
    and McLaughlin, Sean
    and Nguyen, Tat Thang
    and Nguyen, Quang Truong
    and Nipkow, Tobias
    and Obua, Steven
    and Pleso, Joseph
    and Rute, Jason
    and Solovyev, Alexey
    and Ta, Thi Hoai An
    and Tran, Nam Trung
    and Trieu, Thi Diep
    and Urban, Josef
    and Vu, Ky
    and Zumkeller, Roland
  },
  title        = {A Formal Proof of the {Kepler} Conjecture},
  journaltitle = {Forum of Mathematics, Pi},
  volume       = {5},
  date         = {2017},
  eid          = {e2},
  doi          = {10.1017/fmp.2017.1},
}

@inproceedings{HFLL,
  author = {
    Hou (Favonia), Kuen-Bang
    and Finster, Eric
    and Licata, Daniel R.
    and Lumsdaine, Peter LeFanu
  },
  title = {
    A Mechanization of the {Blakers--Massey} Connectivity Theorem
    in Homotopy Type Theory
  },
  editor = {
    Grohe, Martin
    and Koskinen, Eric
    and Shankar, Natarajan
  },
  booktitle = {
    Proceedings of the 31st Annual {ACM/IEEE} Symposium
    on Logic in Computer Science
  },
  series    = {LICS '16},
  date      = {2016},
  pages     = {565--574},
  publisher = {Association for Computing Machinery},
  location  = {New York, NY, USA},
  isbn      = {978-1-4503-4391-6},
  doi       = {10.1145/2933575.2934545},
}

@online{johnson-hopf,
  author  = {Johnson, Niles},
  title = {A visualization of the Hopf fibration},
  urldate = {2026-08-14},
  url = {https://nilesjohnson.net/hopf.html}
}

@inproceedings{Korniowicz:jordan,
  title={A Proof of the Jordan Curve Theorem via the Brouwer Fixed Point Theorem},
  author={Artur Korniłowicz},
  year={2008},
  url={https://api.semanticscholar.org/CorpusID:15229462}
}

@online{lean-eval,
  author        = {{Lean Community}},
  title         = {{Lean AI Formalization Leaderboard}},
  url           = {https://lean-lang.org/eval/},
  urldate       = {2026-08-14}
}

@online{lean-eval-pi3-s2,
  author        = {Morrison, Kim},
  title         = {{$\pi_3$ of the $2$-sphere is $\mathbb{Z}$}},
  organization  = {Lean Focused Research Organization},
  url           = {https://lean-lang.org/eval/problems/pi3_sphere_two_mulEquiv_int/},
  urldate       = {2026-08-14}
}

@inproceedings{LicataShulman2013,
  author       = {Daniel R. Licata and Michael Shulman},
  title        = {Calculating the Fundamental Group of the Circle in Homotopy Type Theory},
  booktitle    = {2013 28th Annual ACM/IEEE Symposium on Logic in Computer Science},
  pages        = {223--232},
  publisher    = {IEEE},
  year         = {2013},
  doi          = {10.1109/LICS.2013.28},
  eprint       = {1301.3443},
  archivePrefix= {arXiv}
}

@online{Licata2013HoTTProgress,
  author  = {Licata, Daniel R.},
  title   = {Homotopy Theory in Type Theory: Progress Report},
  date    = {2013-05-20},
  url     = {https://homotopytypetheory.org/2013/05/20/homotopy-theory-in-type-theory-progress-report/},
  urldate = {2026-08-14},
}

@misc{lurie2003,
      title={On Infinity Topoi}, 
      author={Jacob Lurie},
      year={2003},
      eprint={math/0306109},
      archivePrefix={arXiv},
      primaryClass={math.CT},
      url={https://arxiv.org/abs/math/0306109}, 
}

@inproceedings{mathlib,
  author        = {{The mathlib Community}},
  shortauthor   = {{Mathlib Community}},
  sortname      = {Mathlib Community},
  title         = {The {L}ean {M}athematical {L}ibrary},
  booktitle     = {Proceedings of the 9th {ACM} {SIGPLAN} International Conference on Certified Programs and Proofs},
  series        = {CPP 2020},
  publisher     = {ACM},
  address       = {New Orleans, LA, USA},
  year          = {2020},
  month         = {1},
  pages         = {367--381},
  doi           = {10.1145/3372885.3373824}
}

@book{introhott,
  author    = {Rijke, Egbert},
  title     = {Introduction to Homotopy Type Theory},
  series    = {Cambridge Studies in Advanced Mathematics},
  number    = {219},
  publisher = {Cambridge University Press},
  date      = {2025},
  isbn      = {9781108844161},
  doi       = {10.1017/9781108933568},
}

@software{agda-unimath,
  author        = {Rijke, Egbert and Stenholm, Elisabeth and Prieto-Cubides, Jonathan and Bakke, Fredrik and {others}},
  title         = {{The agda-unimath library}},
  year          = {2021},
  url           = {https://github.com/UniMath/agda-unimath/},
  urldate       = {2026-08-14},
  license       = {MIT},
}

@online{shulman2015,
  author  = {Shulman, Michael},
  title   = {The Univalent Perspective on Classifying Spaces},
  date    = {2015-01-19},
  url     = {https://golem.ph.utexas.edu/category/2015/01/the_univalent_perspective_on_c.html},
  urldate = {2026-08-14},
}

@misc{shulman-infinity-topos,
      title={All $(\infty,1)$-toposes have strict univalent universes}, 
      author={Michael Shulman},
      year={2019},
      eprint={1904.07004},
      archivePrefix={arXiv},
      primaryClass={math.AT},
      url={https://arxiv.org/abs/1904.07004}, 
}

@article{Shulman2024Strange,
  author        = {Shulman, Michael},
  title         = {Strange New Universes: Proof Assistants and Synthetic Foundations},
  journaltitle  = {Bulletin of the American Mathematical Society},
  series        = {New Series},
  volume        = {61},
  number        = {2},
  date          = {2024-04},
  pages         = {257--270},
  doi           = {10.1090/bull/1830}
}

@book{hottbook,
  author        = {{The Univalent Foundations Program}},
  shortauthor   = {{Univalent Foundation Program}},
  sortname      = {Univalent Foundations Program},
  title         = {Homotopy Type Theory: Univalent Foundations of Mathematics},
  publisher     = {Institute for Advanced Study},
  address       = {Princeton, NJ},
  year          = {2013},
  url           = {https://homotopytypetheory.org/book}
}

@software{Voevodsky2011UnivalentBasics,
  author       = {Voevodsky, Vladimir},
  title        = {Univalent Basics},
  subtitle     = {Part of the \texttt{Foundations} Library for Coq},
  date         = {2010-02/2011-09},
  organization = {UniMath},
  url          = {https://github.com/UniMath/Foundations/blob/master/Generalities/uu0.v},
  urldate      = {2026-08-14},
}

@article{Voevodsky,
  author        = {Voevodsky, Vladimir},
  title         = {On the origins and motivations of univalent foundations},
  url           = {https://www.ias.edu/ideas/2014/voevodsky-origins},
  journal       = {The Institute Letter Summer 2014},
  year          = {2014}
}
\endgroup

\end{document}